\documentclass[12pt]{article}
\usepackage{amssymb,amsmath, amsthm, amscd,ifthen}
\usepackage[dvips]{graphicx}
\usepackage{indentfirst}

\renewcommand{\le}{\leqslant}
\renewcommand{\ge}{\geqslant}

\begin{document}

\title{On hypergraph cliques with chromatic number 3\footnote{This work was supported by 
the grant of RFBR N 09-01-00294.}}
\author{D.D. Cherkashin}
\date{}
\maketitle

This work is devoted to a problem in extremal hypergraph theory, which goes back to P. Erd\H{o}s and L. Lov\'asz (see \cite{EL}). Before 
giving an exact statement of the problem, we recall some definitions and introduce some notation. 

Let $ H = (V,E) $ be a hypergraph without multiple edges. We call it {\it $ n $-uniform}, if any of its edges has cardinality $ n $: for every
$ e \in E $, we have $ |e| = n $. By the {\it chromatic number} of a hypergraph $ H = (V,E) $ we mean the minimum number 
$ \chi(H) $ of colors needed to paint 
all the vertices in $ V $ so that any edge $ e \in E $ contains at least two vertices of some different colors. Finally, a hypergraph is said 
to form a {\it clique}, if its edges are pairwise intersecting. 

In 1973 Erd\H{o}s and Lov\'asz noticed that if an $ n $-uniform hypergraph $ H = (V,E) $ forms a clique, then $ \chi(H) \in \{2,3\} $. They also 
observed that in the case of $ \chi(H) = 3 $, one certainly has $ |E| \le n^n $ (see \cite{EL}). Thus, the following definition has been 
motivated: 
$$
M(n) = \max \{|E|: ~ \exists ~ {\rm an} ~ n-{\rm uniform} ~ {\rm clique} ~ H = (V,E) ~ {\rm with} ~ \chi(H) = 3\}.
$$
Obviously such definition has no sense in the case of $ \chi(H) = 2 $. 

\vskip+0.2cm

\noindent{\bf Theorem 1 (P. Erd\H{o}s, L. Lov\'asz, \cite{EL}).} {\it The inequalities hold 
$$
n! \left(\frac{1}{1!} + \frac{1}{2!} + \ldots + \frac{1}{n!}\right) \le M(n) \le n^n.
$$}

\vskip+0.2cm

Almost nothing better has been done during the last 35 years. In the book \cite{TJ} the estimate $ M(n) \le \left(1-\frac{1}{e}\right) n^n $ 
is mentioned as ``to appear''. However, we have not succeeded in finding the corresponding paper. 

At the same time, another quantity $ r(n) $ was introduced in \cite{Lov}:
$$
r(n) = \max \{|E|: ~ \exists ~ {\rm an} ~ n-{\rm uniform} ~ {\rm clique} ~ H = (V,E) ~ {\rm s.t.} ~ \tau(H) = n\},
$$
where $ \tau(H) $ is the {\it covering number} of $ H $, i.e., 
$$
\tau(H) = \min \{|f|: ~ f \subset V, ~ \forall ~ e \in E ~~ f \cap e \neq \emptyset\}.
$$
Clearly, for any $ n $-uniform clique $ H $, we have $ \tau(H) \le n $ (since every edge forms a cover), and if $ \chi(H) = 3 $, then 
$ \tau(H) = n $. Thus, $ M(n) \le r(n) $. Lov\'asz noticed that for $ r(n) $ the same estimates as in Theorem 1 apply and conjectured that 
the lower estimate is best possible. In 1996 P. Frankl, K. Ota, and N. Tokushige (see \cite{FOT}) disproved this conjecture and showed that 
$ r(n) \ge \left(\frac{n}{2}\right)^{n-1} $. 

We discovered a new upper bound for the initial value $ M(n) $. 

\vskip+0.2cm

\noindent{\bf Theorem 2.} {\it There exists a constant $ c > 0 $ such that
$$
M(n) \le c n^{n-\frac{1}{2}} \ln n.
$$}

\vskip+0.2cm

We shall prove Theorem 2 in the next section. 

\section{Proof of Theorem 2}

We shall proceed by citing or proving successive propositions that will eventually lead us to the proof of the theorem. 

\vskip+0.2cm

\noindent{\bf Proposition 1 (P. Erd\H{o}s, L. Lov\'asz, \cite{EL}).} {\it Let $ H = (V,E) $ be an $ n $-uniform clique 
with $ \chi(H) = 3 $. Let $ k $ be an arbitrary integer such that $ 1 \le k \le n $. Take any set $ W \subseteq V $ of cardinality $ k $. 
Let $ E(W) $ denote the set of all edges $ e \in E $ such that $ W \subseteq e $. Then $ |E(W)| \le n^{n-k} $.}

\vskip+0.2cm

Note that in particular, the degree $ {\rm deg}~ v $ of any vertex $ v \in V $ does not exceed $ n^{n-1} $ (here $ k = 1 $). This fact 
entails immediately the estimate $ M(n) \le n^n $. Although we suppose to prove a much better bound, we shall frequently use Proposition 1 
during the proof. 

To any $ n $-uniform hypergraph $ H = (V,E) $ we assign the set 
$$
B(H) = \left\{v \in V: ~ {\rm deg}~v > \frac{|E|}{n^2}\right\}.
$$

\vskip+0.2cm

\noindent{\bf Proposition 2.} {\it Let $ H = (V,E) $ be an $ n $-uniform clique 
with $ \chi(H) = 3 $. Then the two following assertions hold: 
\begin{enumerate}

\item $ |B(H)| < n^3 $; 

\item any edge $ e \in E $ intersects the set $ B(H) $. 

\end{enumerate}
}

\vskip+0.2cm

\paragraph{Proof.} We start by proving the first assertion. Fix an $ H = (V,E) $. Let $ B = B(H) $. We know that 
$ \sum\limits_{v \in V} {\rm deg}~v = n|E| $. Furthermore, 
$$
\sum\limits_{v \in V} {\rm deg}~v \ge \sum\limits_{v \in B} {\rm deg}~v > \frac{|B||E|}{n^2}. 
$$
Thus, $ \frac{|B||E|}{n^2} < n|E| $, which means that actually $ |B| < n^3 $.

To prove the second assertion fix an arbitrary edge $ e $. Since $ H $ is a clique, any $ f \in E $ intersects $ e $. Therefore, 
$ \sum\limits_{v \in e} {\rm deg}~v \ge |E| $. By pigeon-hole principle, there is a vertex $ v \in e $ with $ {\rm deg}~v \ge 
\frac{|E|}{n} \ge \frac{|E|}{n^2} $. So $ v \in B $, and the proof is complete.  

\vskip+0.2cm

\noindent{\bf Proposition 3.} {\it Let $ H = (V,E) $ be an $ n $-uniform clique 
with $ \chi(H) = 3 $. Let $ t \in \{1, \dots, n\} $ and suppose there is an edge $ e \in E $ that 
intersects the set $ B = B(H) $ by at most $ t $ vertices. 
Then there is a vertex $ v \in V $ with $ {\rm deg} ~ v \ge \frac{|E|}{t+1} $.}

\vskip+0.2cm

\paragraph{Proof.} Fix a hypergraph $ H = (V,E) $ and an $ e \in E $ with $ |e \cap B| \le t $. Put $ f = e \cap B $ and 
$ a = |f| \le t $. We know that for any $ v \in f $, one has $ {\rm deg} ~ v > \frac{|E|}{n^2} $. We also
know that for any $ v \in (e \setminus f) $, one has $ {\rm deg} ~ v  \le \frac{|E|}{n^2} $. Finally, we know that $ H $ is a 
clique. Consequently, 
$$
\sum\limits_{v \in f} {\rm deg}~v = \sum\limits_{v \in e} {\rm deg}~v - \sum\limits_{v \in (e \setminus f)} {\rm deg}~v \ge
|E| - (n-a) \frac{|E|}{n^2}.
$$
By pigeon-hole principle, there is a vertex $ v \in f $ with 
$$
{\rm deg}~v \ge \frac{|E| - (n-a) \frac{|E|}{n^2}}{a}.
$$
The right-hand side of the above inequality decreases in $ a \le t $, so that anyway 
$$
{\rm deg}~v \ge \frac{|E| - (n-t) \frac{|E|}{n^2}}{t} = |E| \cdot \frac{n^2-n+t}{n^2t} \ge \frac{|E|}{t+1},
$$
where the last inequality is true, since $ t \in \{1, \dots, n\} $. Proposition 3 is proved. 

\vskip+0.2cm

\noindent{\bf Proposition 4.} {\it Let $ H = (V,E) $ be an $ n $-uniform clique 
with $ \chi(H) = 3 $. Let $ t \in \{1, \dots, n\} $. Then either $ |E| \le t n^{n-1} $, or for any $ e \in E $, we have 
$ |e \cap B(H)| \ge t $.}

\vskip+0.2cm

\paragraph{Proof.} Fix an $ H = (V,E) $ with $ B(H) = B $. Assume that $ |E| > t n^{n-1} $ and there exists an $ e \in E $ such that 
$ |e \cap B| \le t-1 $. By Proposition 3 we can find a vertex $ v $ with $ {\rm deg} ~ v \ge \frac{|E|}{t} > n^{n-1} $, which is in conflict 
with Proposition 1. Thus, our assumption is false, and the proof is complete. 

\vskip+0.2cm

\noindent{\bf Proposition 5.} {\it Let $ H = (V,E) $ be an $ n $-uniform clique 
with $ \chi(H) = 3 $ and $ |E| > n^{n-\frac{1}{2}} $. Suppose that $ n \ge 100 $. 
Then there exist edges $ e,f \in E $ such that $ \left[\sqrt{n}\right] \le |e \cap f| \le n - \left[\sqrt{n}\right] $.} 

\vskip+0.2cm

\paragraph{Proof.} Fix an $ H = (V,E) $ with $ B(H) = B $. Put $ t = \left[\sqrt{n}\right] $. Since $ tn^{n-1} \le n^{n-\frac{1}{2}} < |E| $, 
Proposition 4 tells us that for any $ e \in E $, we have $ |e \cap B| \ge t $.  

Consider the family $ {\cal B}_t = C_B^t $ consisting of all the $ t $-element subsets of the set $ B $. By the first assertion of 
Proposition 2, we have 
$$ 
|{\cal B}_t| = C_{|B|}^t \le |B|^t < n^{3t}.
$$
Also we know that any $ e \in E $ must contain a set $ T \in {\cal B}_t $, since $ |e \cap B| \ge t $.  

At the same time, $ |E| > n^{n-\frac{1}{2}} > n^{5t} $ as $ n \ge 100 $. Thus, taking into account the notation $ E(W) $ from the statement of Proposition 1, 
we see that there exists a set $ T \in {\cal B}_t $ such that $ |E(T)| > n^{5t}/n^{3t} = n^{2t} $.  

Clearly for any $ e,f \in E(T) $, we have $ T \subseteq (e \cap f) $, so that $ |e \cap f| \ge t = \left[\sqrt{n}\right] $. If there exist 
$ e, f \in E(T) $ with $ |e \cap f| \le n - \left[\sqrt{n}\right] = n-t $, then the proposition is proved. Otherwise, every two edges 
from $ E(T) $ intersect by at least $ n-t+1 $ vertices. 

Take an arbitrary $ A \in E(T) $. Put $ s = n-t+1 $ and consider the family $ {\cal A}_s = C_{A}^s $
consisting of all the $ s $-element subsets of the set $ A $. We know that simultaneously {\bf a)} $ |E(T)| > n^{2t} $; {\bf b)} any $ e \in E(T) $
contains a set $ S \in {\cal A}_s $ (since $ |e \cap A| \ge s $); {\bf c)} $ |{\cal A}_s| = C_{|A|}^s = C_n^s = C_n^{t-1} < n^t $. Therefore, 
there is a set $ S \in {\cal A}_s $ such that $ |E(S)| > n^t $. Since $ |S| = s $, by Proposition 1, we have $ |E(S)| \le n^{n-s} = n^{t-1} $, 
which is a contradiction. Proposition 5 is proved. 

\vskip+0.2cm

\noindent{\bf Remark 1.} Note that the proof of Proposition 5 can be easily extended to support the following assertion: 
{\it Let $ H = (V,E) $ be an $ n $-uniform clique 
with $ \chi(H) = 3 $. Let an $ F \subseteq E $ be such that $ |F| > n^{n-\frac{1}{2}} $. Suppose that $ n \ge 100 $. 
Then there exist edges $ f_1,f_2 \in F $ such that $ \left[\sqrt{n}\right] \le |f_1 \cap f_2| \le n - \left[\sqrt{n}\right] $.} 
Note also that a hypergraph $ H' = (V,F) $ does not necessarily have chromatic number 3. It can be bipartite as well. 

\vskip+0.2cm

\noindent{\bf Proposition 6.} {\it Let $ n \in {\mathbb N} $, $ t \in \{1, \dots, n\} $, 
$$
t' = \min \left\{t, 4 \sqrt{n} \ln n\right\}, ~~~
N(t) = (t+1) \left(n - \frac{\sqrt{n}}{4}\right)^{t'-1} n^{n-t'}.
$$
Then $ N(t) = O\left(n^{n-\frac{1}{2}} \ln n\right) $.}  

\vskip+0.2cm

\paragraph{Proof.} First, assume that $ t \le 4 \sqrt{n} \ln n $. Then
$$
N(t) \le (t+1) \cdot n^{t-1} \cdot n^{n-t} = (t+1) n^{n-1} = O\left(n^{n-\frac{1}{2}} \ln n\right), 
$$
and we are done. Now, assume that $ t > 4 \sqrt{n} \ln n $. In this case, 
$$
N(t) \le (n+1) \left(n - \frac{\sqrt{n}}{4}\right)^{t'-1} n^{n-t'} = (n+1) \cdot n^{n-1} 
\left(1 - \frac{1}{4 \sqrt{n}}\right)^{4 \sqrt{n} \ln n-1} = O\left(n^{n-1}\right), 
$$
and we are done again. Proposition 6 is proved. 

\vskip+0.2cm

\noindent{\bf Remark 2.} Note that we may write, say, $ N(t) \le 10 n^{n-\frac{1}{2}} \ln n $ for $ n \ge n_0 $ and all $ t $. 

\vskip+0.2cm

\paragraph{Completion of the proof of the theorem.} Fix an $ n $-uniform clique $ H = (V,E) $ with $ \chi(H) = 3 $ and $ n \ge 
\max \{n_0, 10000\} $. We shall prove that $ |E| \le 10 n^{n-\frac{1}{2}} \ln n $. This will be enough to complete the proof of Theorem 2.

Let 
$$
T = \max \{t: ~ \forall ~ e \in E ~~ |e \cap B| \ge t\}.
$$
By the second assertion of Proposition 2, $ T \in \{1, \dots, n\} $. 

Define $ T' $ in the same way as $ t' $ was defined by $ t $ in Proposition 6. Since $ n \ge 10000 $, we have 
$ T' < n $, and thus $ T' \in \{1, \dots, n-1\} $. Also, since $ n \ge n_0 $, we have $ N(T) \le 
10 n^{n-\frac{1}{2}} \ln n $ (see Remark 2). 

Assume that $ |E| > 10 n^{n-\frac{1}{2}} \ln n $. So we automatically assume that $ |E| > N(T) $. By the definition of the value $ T $,  
there exists an edge $ e \in E $ that intersects $ B $ by at most $ T $ vertices. Hence, by Proposition 3 there is a vertex $ v \in V $ 
with 
$$ 
{\rm deg}~v \ge \frac{|E|}{T+1} > \frac{N(T)}{T+1} = \left(n - \frac{\sqrt{n}}{4}\right)^{T'-1} n^{n-T'}.
$$

Put $ I = \{v\} $, $ i = 1 $. Then 
$$
|E(I)| = {\rm deg}~v > \left(n - \frac{\sqrt{n}}{4}\right)^{T'-1} n^{n-T'} = 
\left(n - \frac{\sqrt{n}}{4}\right)^{T'-i} n^{n-T'}.
\eqno{(1)}
$$

If $ T' = 1 $, then inequality (1) contradicts Proposition 1. Therefore, assume that $ T' \in \{2, \dots, n-1\} $. 

Let inequality (1) serve as the base for an inductive procedure with $ \le T' $ steps. So assume that we have already found a set $ I \subset V $ with 
$ |I| = i \in \{1, \dots, T'-1\} $ (do not forget that $ T' \ge 2 $) and 
$$
|E(I)| > \left(n - \frac{\sqrt{n}}{4}\right)^{T'-i} n^{n-T'}.
$$
We shall prove that either we can take an $ a \in (V \setminus I) $ such that 
$$
|E(I \cup \{a\})| > \left(n - \frac{\sqrt{n}}{4}\right)^{T'-i-1} n^{n-T'},
\eqno{(2)} 
$$
or we can take $ a, b \in (V \setminus I) $ such that 
$$
|E(I \cup \{a,b\})| > \left(n - \frac{\sqrt{n}}{4}\right)^{T'-i-2} n^{n-T'}.
\eqno{(3)}
$$
(Here for all $ i $, $ T' - i - 2 \ge -1 $ and $ i+2 \le T'+1 \le n $, so that the choice of the parameters is correct.) 

Indeed, to prove (2) or (3) let 
$$
E_I = \{e \in E: ~ e \cap I \neq \emptyset\}. 
$$
By Proposition 1 $ |E_I| \le i n^{n-1} $. Hence, $ |E_I| < T' n^{n-1} \le 4 n^{n-\frac{1}{2}} \ln n $. Putting 
$$
E^I = \{e \in E: ~ e \cap I = \emptyset\}
$$
we immediately get the estimate 
$$
|E^I| = |E| - |E_I| > 10 n^{n-\frac{1}{2}} \ln n - 4 n^{n-\frac{1}{2}} \ln n = 6 n^{n-\frac{1}{2}} \ln n > n^{n-\frac{1}{2}}.
$$
Since $ n > 100 $, Remark 1 tells us that there exist $ f_1, f_2 \in E^I $ with 
$ \left[\sqrt{n}\right] \le |f_1 \cap f_2| \le n - \left[\sqrt{n}\right] $.

Since $ H $ is a clique, any edge $ e $ from $ E(I) $ intersects both $ f_1 $ and $ f_2 $. So either $ e $ goes through a vertex  
$ v \in (f_1 \cap f_2) $ or it contains a vertex $ v_1 \in (f_1 \setminus f_2) $ and a vertex $ v_2 \in (f_2 \setminus f_1) $. Formally, 
we may write down the equality 
$$
E(I) = \left(\bigcup_{v \in (f_1 \cap f_2)} E_v(I)\right) \bigcup \left(\bigcup_{v_1 \in (f_1 \setminus f_2)} \bigcup_{v_2 \in (f_2 \setminus f_1)}
E_{v_1,v_2}(I)\right), 
$$
where
$$
E_v(I) = \{e \in E(I): ~ v \in e\}, ~~~
E_{v_1,v_2}(I) = \{e \in E(I): ~ v_1, v_2 \in e\}.
$$

Of course
$$
|E(I)| \le \sum_{v \in (f_1 \cap f_2)} \left|E_v(I)\right| + \sum_{v_1 \in (f_1 \setminus f_2)} \sum_{v_2 \in (f_2 \setminus f_1)}
\left|E_{v_1,v_2}(I)\right|. 
$$

If a summand in the first sum is greater than $ \left(n - \frac{\sqrt{n}}{4}\right)^{T'-i-1} n^{n-T'} $, then (2) is shown: indeed, the 
corresponding $ v $ is contained in this many edges $ e \in E(I) $ that already contain $ I $. If a summand in the second sum is 
greater than $ \left(n - \frac{\sqrt{n}}{4}\right)^{T'-i-2} n^{n-T'} $, then (3) is shown in turn. So suppose that there are no such 
summands. In this case, putting $ k = |f_1 \cap f_2| $ we have
$$
|E(I)| \le k \left(n - \frac{\sqrt{n}}{4}\right)^{T'-i-1} n^{n-T'} + (n-k)^2 \left(n - \frac{\sqrt{n}}{4}\right)^{T'-i-2} n^{n-T'}.
$$

On the other hand, $ |E(I)| > \left(n - \frac{\sqrt{n}}{4}\right)^{T'-i} n^{n-T'} $. Having proved that 
$$
\left(n - \frac{\sqrt{n}}{4}\right)^{T'-i} n^{n-T'} > 
k \left(n - \frac{\sqrt{n}}{4}\right)^{T'-i-1} n^{n-T'} + (n-k)^2 \left(n - \frac{\sqrt{n}}{4}\right)^{T'-i-2} n^{n-T'}   
$$
for any $ k \in \left[\left[\sqrt{n}\right], n - \left[\sqrt{n}\right]\right] $, we would get a contradiction which would complete the 
proof of (2) or (3). 

The needed inequality is equivalent to 
$$
\left(n - \frac{\sqrt{n}}{4}\right)^2 > k \left(n - \frac{\sqrt{n}}{4}\right) + (n-k)^2,
$$
which can be proved by standard analytic calculations. 

Thus, (2) or (3) take place. So after $ \le T' $ steps of the inductive procedure, we get either a set $ I $ of cardinality $ T' $ such that 
$ |E(I)| > n^{n-T'} $ or a set $ I $ of cardinality $ T'+1 $ such that $ |E(I)| > \left(n - \frac{\sqrt{n}}{4}\right)^{-1} n^{n-T'} > 
n^{n-T'-1} $. Both estimates are in conflict with Proposition 1. Consequently, our initial assumption $ |E| > 10 n^{n-\frac{1}{2}} \ln n $
is false, and Theorem 2 is proved.

\end{document}